\documentclass[centertags,12pt,a4paper]{amsart}

\usepackage{latexsym}
\usepackage{amssymb,amsmath}
\usepackage{amscd}

\makeatletter
\@addtoreset{equation}{section}
\makeatother

\newtheorem{theorem}{Theorem}[section]
\newtheorem{proposition}[theorem]{Proposition}
\newtheorem{deff}[theorem]{Definition}
\newtheorem{lemma}[theorem]{Lemma}

\newtheorem*{claim*}{Claim}
\theoremstyle{definition}

\newtheorem*{question*}{Question}

\def\N{\mathbf{N}}

\def\bfq{{\bar{\mathbb F}}_q}

\def\fqc{\mathbb{F}_{q^3}}
\def\fq{\mathbb{F}_q}

\author[F. Pasticci]{Fabio Pasticci}

\thanks{This research was performed within the activity of GNSAGA of the
Italian INDAM, with the financial support of the Italian Ministry MIUR,
project ``Strutture geometriche, combinatorica e loro applicazioni'',
PRIN 2006-2007}

\title[On the genus of a cyclic plane curve over a finite field]{On the genus of a cyclic plane curve over a finite field}

\address{Dipartimento di Matematica e Informatica Universit\`a degli Studi di Perugia, 06123
Perugia, Italy}\email{pasticci@dipmat.unipg.it}

\begin{document}

\begin{abstract}
Cyclic curves, i.e. curves fixed by a cyclic collineation group,
play a central role in the investigation of cyclic arcs in
Desarguesian projective planes. In this paper, the genus of a
cyclic curve arising from a cyclic $k$-arc of Singer type is
computed.
\end{abstract}

\maketitle


\section{Introduction}

Let $PG(2,q)$ be the projective plane over the finite field ${\mathbb F}_q$, $q=p^h$ for some prime $p$ and some $h \in \N$.
A $k-$arc in $PG(2,q)$, $k \geq 3$ is a set of $k$ points, every 3 of which are not collinear.
A $k-$arc in $PG(2,q)$ is said to be complete if and only if it cannot be extended to a $(k+1)-$arc by a point of $PG(2,q)$.
A $k$-arc is called {\em cyclic } \cite{GIU} if it consists of the
points of a point orbit under a cyclic collineation group $G$ of
$PG(2,q)$.
A cyclic $k-$arc is said to be {\em of Singer type} if it consists of
a point orbit under a subgroup of a cyclic Singer group of $PG(2,q)$.

An essential tool in the investigation of $k-$arcs is the following result due to B. Segre (see \cite{HIR}).

 If $q$ is odd, then there exists
a plane curve $\Gamma '$ in the dual plane of $PG(2,q)$ such that:

\begin{enumerate}

\item $\Gamma '$ is defined over ${\mathbb F_q}$.

\item The degree of $\Gamma '$ is $2t$, with $t=q-k+2$ being the number of
1-secants through a point of $K$.

\item The $kt$ 1-secants of $K$ belong to $\Gamma '$.

\item Each 1-secant $\ell$ of $K$ through a point $P\in K$ is counted
twice in the intersection of $\Gamma '$ with $\ell_P$,  where $\ell_P$ denotes the line corresponding to $P$ in the dual plane.

\item The curve $\Gamma '$ contains no 2-secant of $K$.

\item The irreducible components of $\Gamma '$ have multiplicity at most
2, and $\Gamma '$ has at least one component of multiplicity 1.

\item If $k> \frac{2}{3}(q+2)$, then there exists a unique curve in the
dual plane of $PG(2, {\bar {\mathbb F_q}})$ satisfying properties (2), (3), (4), (5).

\end{enumerate}

The investigation of the algebraic envelope of a cyclic $k-$arc of Singer type was initiated by Cossidente and Korchm\'aros in 1998 \cite{CK}, and continued by Giulietti in \cite{GIU}. In \cite{GIU} the terminology of {\em cyclic curve} is introduced to denote the algebraic plane curve defined over $\fq$  corresponding to the algebraic envelope of a cyclic $k-$arc of Singer type.

In this paper we deal with the problem of computing the genus of a cyclic curve. Our main result is the following:

\begin{theorem}
Let ${\mathcal C}$ be a cyclic curve of order $n=2(q-k+2)=2t$. If
$k \geq q-\sqrt{\frac{2}{3}q+\frac{1}{4}}+\frac{9}{4}$, $k^2-k+1
\neq 0$ mod $p$, and ${\mathcal C}$ is irreducible of genus $g$,
then ${\mathcal C}$ is singular and either $g=2 { (t-1)(t-2) }$ or
$g= \frac{(t-1)(t-2)}{2} $ holds.

\end{theorem}

\section{Preliminaries on Singer cycles}

Let $ \fqc  $  be a cubic extension of $\fq$. Following Singer \cite{SIN}, we identify the projective plane  $PG(2,q)$ with $\fqc$ $\mod$ $\fq$. This means that points of $PG(2,q)$ are non-zero elements  of $\fqc$ and two elements $x,y \in \fqc$  represent the same point of $PG(2,q)$ if and only if $x/y \in \fq$.
Let $\omega$ be a primitive element of $\fqc$
and $$p(x)=x^3-ax^2-bx-c$$ its minimal polynomial over $\fq$.

The matrix

$$
C=
\left(
\begin{array}{ccc}
0&1&0 \\
0&0&1 \\
c&b&a \\
\end{array}
\right)
$$

induces a linear collineation $\phi$ of $PG(2,q)$ of order $q^2+q+1$ called a Singer cycle of $PGL(3,\fq )$.
All Singer cycles of $PGL(3, \fq)$ form a single conjugacy class and the matrix $C$ is conjugate in $GL(3,\fq)$ to the diagonal matrix

$$
D=
\left(
\begin{array}{ccc}
\omega &  0        &   0  \\
\omega & \omega^q  &   0  \\
0      &   0       & \omega^{q^2} \\
\end{array}
\right)
$$

by the matrix

$$
E=
\left(
\begin{array}{ccc}
   1      &  1          &   1  \\
\omega    & \omega^q    & \omega^{q^2}  \\
\omega^2  & \omega^{2q} & \omega^{2q^2} \\
\end{array}
\right)
$$

Let $\sigma$ denote the linear collineation of $PG(2,q^3)$  induced by $D$. It fixes the points $E_0=(1,0,0)$, $E_1=(0,1,0)$ and $E_2=(0,0,1)$.
The linear collineation $T$ of $PG(2,q^3)$ defined by

$$(X_0,X_1,X_2) \to (X_2,X_0,X_1)$$

has order 3 and acts on the points $E_0, E_1, E_2$ as the cycle $(E_0\; E_1\; E_2)$

\begin{proposition}\label{conj}
Cyclic Singer groups of $PG(2,q)$ are equivalent under conjugation by the
elements of $PGL(3,q)$.
\end{proposition}
\begin{proof}
See  \cite{SZ} and \cite{HIR}, Corollary 4 to Theorem. 4.2.1.
\end{proof}

\begin{deff}
A {\em $k$-arc of Singer type} in $PG(2,q)$ is a $k$-arc which consists of
a point orbit under a subgroup of a cyclic Singer group of $PG(2,q)$.
\end{deff}

From now on $K$ will be a $k$-arc of Singer type in $PG(2,q)$;
obviously, $k$ divides $q^2+q+1$. Since by Proposition \ref{conj} two orbits
of the same size under subgroups of cyclic Singer groups are
projectively equivalent, we will assume without loss of generality that
\begin{equation}\label{arc}
K=\{1, \omega^{\frac{q^{2}+q+1}{k}}, \ldots ,
\omega^{(k-1)\frac{q^{2}+q+1}{k} } \}.
\end{equation}

Let $K$ be a $k-$arc of $PG(2,q)$, we consider the unique envelope $\Gamma '$ of $K$ as a plane algebraic curve in $PG(2, \bfq)$ defined over $\fq$ and we will denote it by $\Gamma_{2t}$.

Following \cite{GIU}, we will study a curve  ${\mathcal C}_{2t}$ projectively
equivalent to $\Gamma_{2t}$ in $PG(2,q^{3})$.
Let $\varphi=[L]$ be the element of $PGL(3,q^{3})$ such that
\[
LBL^{-1}= \left(
              \begin{array}{ccc}
           \omega^q & 0 & 0 \\
           0 & \omega^{q^2} & 0 \\
           0 & 0 & \omega
              \end{array}
          \right)
\]
and $\varphi(1,0,0)=(1,1,1)$. Let  $\Pi$ be the image of $PG(2,q)$ by
$\varphi$; $\Pi$ is a subplane of
$PG(2,q^{3})$ isomorphic to $PG(2,q)$.
We can consider $\alpha := \varphi \sigma \varphi^{-1}$,
collineation of $\Pi$; we have that if $a:=\omega^{q-1}$, then
\[ \alpha :
 \left\{
    \begin{array}{cl}
    \rho x_{1}^{'}  = & ax_{1} \\
    \rho x_{2}^{'} = & a^{q+1}x_{2} \\
    \rho x_{3}^{'} = & x_{3}
    \end{array}
    \right. ,\:\:\:\:\:\:\:\:\: \rho \in \fqc^* .
\]
So the points of $\Pi$ are those in
the orbit of $(1,1,1)$ under the action of
$<\alpha>$:
\[
\Pi=\{ (a^{\imath},a^{\imath(q+1)},1) \mid \imath=0,1, \ldots , q^{2}+q \}.
\]
By \ref{arc} $K$ is fixed by
$\sigma^{\frac{q^{2}+q+1}{k}}$
and $\tau$, where $\tau$ is defined by $\tau(\omega^{\imath}):=
\omega^{q\imath}$.
We are interested in $\beta := \varphi \sigma^{\frac{q^{2}+q+1}{k}}
\varphi^{-1}$
and $\delta := \varphi \tau \varphi^{-1}$, collineations of $\Pi$. We have that if
$b:=a^{\frac{q^{2}+q+1}{k}}$  then
\[
\beta :
 \left\{
    \begin{array}{cl}
    \rho x_{1}^{'}  = & bx_{1} \\
    \rho x_{2}^{'} = & b^{q+1}x_{2} \\
    \rho x_{3}^{'} = & x_{3}
    \end{array}
    \right. ,\:\:\:\:\:\:\:\rho \in \fqc^*
\]
and
\[
 \delta :
 \left\{
    \begin{array}{cc}
    \rho x_{1}^{'}  = & x_{1}^{q} \\
    \rho x_{2}^{'} = & x_{2}^{q} \\
    \rho x_{3}^{'} = & x_{3}^{q}
    \end{array}
    \right. , \:\:\:\:\:\:\: \rho \in \fqc^*.
\]
$\delta$ is a non-linear collineation
of $PG(2,q^{3})$, but it acts on  $\Pi$ as
$\mu^{-1}$, where $\mu $ is defined by
\[ \mu :
 \left\{
    \begin{array}{cc}
    \rho x_{1}^{'}  = & x_{3} \\
    \rho x_{2}^{'} = & x_{1} \\
    \rho x_{3}^{'} = & x_{2}
    \end{array}
    \right. ,\:\:\:\:\:\:\: \rho \in \fqc^*.
\]

In \cite{GIU} the following proposition is proved.
\begin{proposition}\label{properties} 
Let $\Gamma_{2t}$ be the algebraic curve associated to the envelope of $K$
and let  ${\mathcal C}_{2t}$ be its image by  $\varphi$. Then
${\mathcal C}_{2t}$ has the following properties:
\begin{enumerate}
 \item it is preserved by $\beta$;
 \item it is preserved by $\mu$;
 \item it has degree $2(q-k+2)$;
 \item it has no fundamental line as a component;
 \item it is defined over ${\mathbb F_q}$.
\end{enumerate}

\end{proposition}

\begin{deff}
Let $\omega$ be a primitive element of $\fqc$, $k$ be a divisor of
$q^2+q+1$ such that $k>\frac{2}{3}(q+2)$, and suppose $\beta$ and $\mu$
defined as above. An algebraic plane curve
defined over ${\mathbb F_q}$ is called {\em cyclic} if it satisfies {\rm (1), (2),
(3), (4)} of Proposition\rm{ \ref{properties}}.
\end{deff}

For the rest of the section ${\mathcal C}$ will denote a cyclic curve of
degree $n=2(q-k+2)$.
\begin{proposition}\label{CK}
Each vertex of the fundamental triangle is a 2-fold cuspidal point of
${\mathcal C}$ such that one of the fundamental line through  the vertex is
the tangent and has intersection multiplicity $n-2$ with ${\mathcal C}$ at
the vertex, namely,
letting  $A_{1}:=(1,0,0)$, $A_{2}:=(0,1,0)$, $A_{3}:=(0,0,1)$,
\[
I(A_2;{\mathcal C}\cap \{x_{1}=0\})=
I(A_3;{\mathcal C}\cap \{x_{2}=0\})=
I(A_1;{\mathcal C}\cap \{x_{3}=0\})=n-2.
\]
Besides, the only two possibilities for branches of ${\mathcal C}$ centered
at $A_{\imath}$ ($\imath =0,1,2$) are the following:
\begin{enumerate}
\item there exist two linear branches (not necessarily distinct)
centered at $A_{\imath}$ and the tangent meets each of them with multiplicity
$\frac{n}{2}-1$;
\item there exists a unique quadratic branch centered at $A_{\imath}$.
\end{enumerate}
\end{proposition}
\begin{proof}
See \cite{CK}, Prop. 5.
\end{proof}
Note that $<\mu> $ preserves ${\mathcal C}$ and acts transitively
on the vertices of the fundamental triangle; so the number of
branches through $A_\imath$ and their characteristics do not
depend on $\imath$. According to \cite{GIU}, if (1) of the above
Proposition holds, then  {$\mathcal C$} is said to be cyclic {\em
of the first type}, otherwise {\em of the second type}.

\section{On the Genus of a Cyclic Curve}
In this section we will prove that if
$k \geq q-\sqrt{\frac{2}{3}q+\frac{1}{4}}+\frac{9}{4}$ and if
the envelope of the $k$-arc of Singer type $K$ in $PG(2,q)$ is irreducible,
than we can establish its genus. More precisely the following theorem
holds (the notation will be as in Section 2).
\begin{theorem}\label{gen}
Let ${\mathcal C}$ be a cyclic curve of order $n=2(q-k+2)=2t$, as defined
in
Section 2.  If $k \geq q-\sqrt{\frac{2}{3}q+\frac{1}{4}}+\frac{9}{4}$,
$k^2-k+1 \neq 0$ mod $p$, and
${\mathcal C}$ is irreducible of genus $g$, then
\begin{itemize}
\item if ${\mathcal C}$ is of the first type
then $2g-2=4t^{2}-12t+6$;
\item if ${\mathcal C}$ is of the second type then $2g-2=t^{2}-3t$.
\end{itemize}
\end{theorem}

Through the rest of the section we will prove Theorem \ref{gen}.
Let  ${\mathcal C}$ be a cyclic curve
and let $g(x,y)=0$ its minimal equation. Let $g(x,y)=\varphi_{\jmath}(x,y)+
\ldots +\varphi_{n}(x,y)$, $\varphi_{u}$ homogeneous of degree $u$; we will
denote the generic term of $\varphi_{u}(x,y)$ by
$a_{\imath,u-\imath}x^{\imath}y^{u-\imath}$.
In \cite{CK} the following four lemmas are proved.
\begin{lemma}
For every $m$, $(0 \leq m \leq n)$, $g(x,y)$ has at most one term of degree
$m$.
\end{lemma}
\begin{lemma}
For every $l$, $(0 \leq l \leq n)$, $g(x,y)$ has at most one term of degree
$l$ in $x$ and at most one term of degree $n-l$ in $y$.
\end{lemma}
\begin{lemma}\label{lkj}
For any two integers $l,\: m$ $(0 \leq l \leq m \leq n)$, we have
\[
a_{l, m-l}=\epsilon a_{n-m,l}=\epsilon^{2}a_{m-l,n-m},
\]
with $\epsilon^{3}=1$.
\end{lemma}
\begin{lemma}\label{ggmm}
Let $a_{l,m-l}x^{l}y^{m-l}$ be a term of $g(x,y)$ different from zero; then
\[
m \equiv (s-1)l+2 \:\:\:\:mod \:\:k.
\]
\end{lemma}
Now we can prove the following proposition.
\begin{proposition} Let ${\mathcal C}$ be a cyclic curve of order
$n=2(q-k+2)=2t$. If $k \geq q-\sqrt{\frac{2}{3}q+\frac{1}{4}}+\frac{9}{4}$
then ${\mathcal C}$ has equation
\[
g(x,y)=y^{2}+\epsilon_{1}x^{2}y^{2t-2}+\epsilon_{1}^{2}x^{2t-2}+c
(x^{t}y^{t-1}+\epsilon_{2}x^{t-1}y+\epsilon_{2}^{2}xy^{t})=0,
\]
with $\epsilon_{1}^{3}=\epsilon_{2}^{3}=1$.
\end{proposition}
\begin{proof} For any term $a_{l,m-l}x^{l}y^{m-l}$ of $g(x,y)$ there exist
by Lemma \ref{lkj} two other terms of $g(x,y)$ of type
$a_{n-m,l}x^{n-m}y^{l}$ and $a_{m-l,n-m}x^{m-l}y^{n-m}$. By permutating
indexes we may assume that $0 \leq l \leq \frac{1}{3}n$; so
\[
(t-1)l+2 \leq \frac{2}{3}(t-1)t+2;
\]
the hypothesis concerning $k$ yields $\frac{2}{3}(t-1)t+2 \leq k$. By Lemma
\ref{ggmm} we have $m=(t-1)l+2$ and since $m \leq 2t$ we have $l \leq 2$.
So the couple $(l,m-l)$ is equal to $(0,2)$, $(1,t)$ or $(2,2t-2)$ and the
proposition is proved (use  Lemma \ref{lkj} again).
\end{proof}
We will compute the genus of ${\mathcal C}$ by applying the famous
Hurwitz Theorem. Let ${\mathcal L}:={\bar {\mathbb F}_q}$ and let
$\Sigma={\mathcal L}({\mathcal C})$ the field of rational functions of
${\mathcal C}$. Let $x$ and $y$ be elements of $\Sigma$ such that
$\Sigma={\mathcal L}(x,y)$. Following Seidenberg's book approach to algebraic
curves (see \cite{SEI}) we define a rational transformation $\Phi$
of ${\mathcal C}$ in $PG(2,{\mathcal L})$ by choosing two elements $x'$ and
$y'$ in $\Sigma$:
\[
x':=\frac{y}{x^{t-1}}, \:\:\:\:\:\:\:\:\: y':=\frac{y^{t-1}}{x^{t-2}}.
\]
The image ${\mathcal C}'$ of ${\mathcal C}$ by $\Phi$ is a plane algebraic
curve of genus 0; for  $x'$ and $y'$ satisfy
\[
{x'}^{2}+\epsilon_{1}{y'}^{2}+\epsilon_{1}^{2}+c{y'}+c\epsilon_{2}{x'}+
c\epsilon_{2}^{2} {x'}{y'}=0,
\]
so ${\mathcal C}'$ is a conic or a line. Let $\Sigma '={\mathcal L}
(x',y')$.
\begin{lemma}\label{degree}
The degree $[\Sigma : \Sigma ']$ of the extension $\Sigma : \Sigma '$ is
$t^{2}-3t+3$ or $2(t^{2}-3t+3)$.
\end{lemma}
\begin{proof} Note that the only branches of ${\mathcal C}$ whose image
by $\Phi$ is centered at a point of the line $y'=0$ are those centered
at $A_{3}$.
If ${\mathcal C}$ is of the first type, then at $A_{3}$ are centered two
branches $\gamma_1$, $\gamma_2$ of ${\mathcal C}$ with
$\gamma_{1}=(\tau,a_{0}\tau^{\frac{n}{2}-1}+ \ldots)$,
$\gamma_{2}=(\tau,a_{0}'\tau^{\frac{n}{2}-1}+ \ldots)$ and $a_0, a_{0}'\neq 0$;
$\gamma_{1}$ and $\gamma_{2}$ are transformed by $\Phi$ in
$\gamma_{1}'$ and $\gamma_{2}'$, branches with (imprimitive) representations
$((a_{0}+\ldots),\tau^{t^{2}-3t+3}(a_{0}+\ldots)^{t-1})$ and
$((a_{0}'+\ldots),\tau^{t^{2}-3t+3}(a_{0}'+\ldots)^{t-1})$ respectively.
Since $\gamma_{1}'$ and $\gamma_{2}'$ are branches of a conic or a line,
their intersection multiplicity $I$ with the line $y'=0$ has to be 1 or 2;
but $t^{2}-3t+3$ is odd so $I=1$ and the ramification index
of $\gamma_{1}$ and $\gamma_{2}$ with respect to $\Phi$ is $t^{2}-3t+3$.
Therefore $[\Sigma : \Sigma ']=t^{2}-3t+3$ or
$[\Sigma : \Sigma ']=2(t^{2}-3t+3)$ according to whether
$\gamma_{1}'$ and $\gamma_{2}' $ are distinct or not.
If ${\mathcal C}$ is of the second type, then at $A_{3}$ is centered exactly
one branch $\gamma$ of ${\mathcal C}$ with
$\gamma =(\tau^{2}, b_{0}\tau^{n-2}+ \ldots)$, $ b_{0}\neq 0 $.
$\gamma$ is transformed by $\Phi$ in $\gamma '$ whose (imprimitive)
representation is of type $((a_{0}+ \ldots), \tau^{2(t^{2}-3t+3)}(a_{0}+
\ldots)^{t-1})$; the ramification index of $\gamma$ is equal to the degree
of $\Sigma : \Sigma '$ since $\gamma$ is the only branch over $\gamma '$;
it can be $2(t^{2}-3t+3)$ or $t^{2}-3t+3$ as before and so we are done.
\end{proof}
We will calculate the order of the different $D$ of $\Sigma : \Sigma '$.
We recall that $D$ is the divisor of $\Sigma$
\[
D:= \sum_{\gamma}D_{\gamma}
\]
where $\gamma$ runs in the set of all branches of ${\mathcal C}$ and
$D_{\gamma}$ is defined as follows: if $(\gamma_{1}(\tau),\gamma_{2}(\tau))$
is a primitive representation of $\gamma$ and if $(\psi_{1}(z),\psi_{2}(z))$
is a primitive representation of $\gamma '$, the image of $\gamma$ by $\Phi$,
with $z=z(\tau)=\tau^{v}+ \ldots$ and
$(\psi_{1}(z(\tau)),\psi_{2}(z(\tau)))=
(\Phi(\gamma_{1}(\tau),\gamma_{2}(\tau)))$, then
\[
D_{\gamma}=ord_{\tau}\:\:\:\frac{dz}{d\tau}.
\]
\begin{lemma}\label{56}
Let $u$ be an element of $\mathcal L$ such that $u^{t^{2}-3t+3}=1$. Then
the linear collineation of $PG(2,{\mathcal L})$ $\eta$ defined by
\[ \eta :
 \left\{
    \begin{array}{cl}
    \rho x_{1}^{'}  = & ux_{1} \\
    \rho x_{2}^{'} = & u^{t-1}x_{2} \\
    \rho x_{3}^{'} = & x_{3}
    \end{array}
    \right. , \:\:\:\:\:\:\: \rho \in {\mathcal L}^*
\]
preserves $\mathcal C$.
\end{lemma}
\begin{proof}
The proof is a simple computation.
\end{proof}
\begin{lemma}
If $k^2-k+1 \neq 0 $ mod $p$ then there exist exactly  $t^{2}-3t+3$ distinct
elements $u$ in $\mathcal L$ such that $u^{t^{2}-3t+3}=1$.
\end{lemma}
\begin{proof}
The polynomial $X^{t^{2}-3t+3}-1 $ is inseparable
in ${\mathcal L}[X]$ since the characteristic $p$ of ${\mathcal L}$ does not
divide $t^{2}-3t+3$; for
$t^{2}-3t+3 = k^{2}-k+1$ {\em mod} $p$.
\end{proof}
\begin{proposition}\label{imp}
If $k^2-k+1 \neq 0$ mod $p$ then  $[\Sigma : \Sigma ']=t^{2}-3t+3$ and the
order of the different $D$ of $\Sigma :\Sigma '$ is
$6(t^{2}-3t+2)$ or $3(t^{2}-3t+2)$ according to whether
${\mathcal C}$ is of the first type or not.
\end{proposition}
\begin{proof}
Let $\gamma$ be a branch of $\mathcal C$ centered at a
point $(x_{0}, y_{0})$ not belonging to any fundamental line
and let $\gamma '$ be its image by $\Phi$; $\gamma '$ is
centered at
$(x_{1},y_{1}):=(\frac{y_{0}}{x_{0}^{t-1}},\frac{y_{0}^{t-1}}{x_{0}^{t-
2}})$.
For every $u$ in $\mathcal L$ such that $u^{t^{2}-3t+3 }=1$ the point
$Q:=(ux_{0},u^{t-1}y_{0})$ is a point of $\mathcal C$ and every branch
centered at $Q$ has as image $\gamma '$, the only branch of ${\mathcal C}'$
centered at $(x_{1},y_{1})$. Moreover, it is easy to see that
all the branches of $\mathcal C$ that lie over $\gamma '$ are centered at
points of type $(ux_{0},u^{t-1}y_{0})$ with $u^{t^{2}-3t+3}=1$. So the set
of branches $\gamma '$ of ${\mathcal C}'$ such that there exist exactly
$t^{2}-3t+3$ distinct branches of $\mathcal C$ that lie over $\gamma '$ is
infinite; for otherwise there would be an infinite set of singular points
of $\mathcal C$. Suppose $[\Sigma : \Sigma ']=2(t^{2}-3t+3)$; then there
exist an infinite set of branches $\gamma$ of $\mathcal C$ such that
$D_{\gamma}$ is greater than 0, a contradiction. So by Lemma \ref{degree}
the first part of the statement is proved. Now it is clear that if
$\gamma$ is a branch of $\mathcal C$ centered at a point not belonging to any
fundamental line then $D_\gamma=0$; on the other hand, if $\gamma$ is
centered at a vertex of the fundamental triangle, we have that (with
notation as before Lemma \ref{56}) $z(\tau)=\tau^{t^{2}-3t+3}+ \ldots$ so
\[
\frac{dz}{d\tau}=(t^{2}-3t+3)\tau^{t^{2}-3t+2}+ \ldots;
\]
since $p$ does not divide $t^{2}-3t+3$ we have
$ord \frac{dz}{d\tau}=t^{2}-3t+2$ and we are done.
\end{proof}
We recall the statement of the Hurwitz Theorem.
\begin{theorem}
Let $\mathcal G$ be an irreducible algebraic curve over the algebraically
closed field  ${\mathcal L}$ and let ${\mathcal G}'$ be the image of
$\mathcal G$ by a rational transformation. Let $\Sigma:={\mathcal
L}({\mathcal G})$ and $\Sigma ' :={\mathcal L}({\mathcal G}')$, $\Sigma '
\subseteq \Sigma$. Let $D$ be the different of $\Sigma : \Sigma '$, $g$ be
the genus of $\mathcal G$  and $g'$ be the genus of ${\mathcal G}'$. Then
\begin{equation}\label{Hur}
2g-2=[\Sigma : \Sigma '](2g'-2)+ ord(D).
\end{equation}
\end{theorem}
Finally we prove the main result of the section.
\begin{proof}
{\em (Theorem \ref{gen})} We apply (\ref{Hur}) and Proposition \ref{imp}; so
we have that if $\mathcal C$ is of the first type then
\[
2g-2=(t^{2}-3t+3)(-2)+6(t^{2}-3t+2)=4t^{2}-12t+6,
\]
otherwise
\[
2g-2=(t^{2}-3t+3)(-2)+3(t^{2}-3t+2)=t^{2}-3t.
\]
\end{proof}

\end{document}